\documentclass[a4paper,11pt,twoside]{article}
\usepackage{pst-fill,pst-grad}
\usepackage{textcomp}
\usepackage[english]{babel}
\usepackage{graphicx}
\usepackage{amsmath}
\usepackage{float}
\usepackage[matrix,arrow,curve]{xy}
\usepackage{pstricks}
\usepackage{mathtools}
\usepackage{amsmath,amsfonts,verbatim,afterpage,theorem,euscript,mathrsfs,amssymb}
\usepackage{amsfonts}
\usepackage{amssymb}
\usepackage{array}
\usepackage{dsfont}
\usepackage{hyperref}
\newtheorem{Theoreme}{Theorem}

\newtheorem{Lemme}{Lemma}[section]

\numberwithin{equation}{section}
\def\vu{\vec{u}}

\def\vf{\vec{f}}
\def\vg{\vec{g}}

\def\vn{\vec{\nabla}}

\def\Rt{\mathbb{R}^3}
\def\R{\mathbb{R}}

\title{\bf   {L}iouville type theorems for stationary
{N}avier-{S}tokes equations with Lebesgue spaces of variable exponent }
\usepackage{authblk}
\author{Diego Chamorro\footnote{\emph{diego.chamorro@univ-evry.fr} (corresponding author)} }
\author{Gast\'on Vergara-Hermosilla\footnote{\emph{gaston.vergarahermosilla@univ-evry.fr}} }
\affil{\footnotesize LaMME, Univ. Evry, CNRS, Universit\'e Paris-Saclay, 91025, Evry, France.}
\begin{document}
\maketitle
\begin{scriptsize}
\abstract{In this article we study some Liouville-type theorems for the stationary 3D Navier-Stokes equations. These results are related to the uniqueness of weak solutions for this system under some additional information over the velocity field, which is usually stated in the literature in terms of Lebesgue, Morrey or $BMO^{-1}$ spaces. Here we will consider Lebesgue spaces of variable exponent which will provide us with some interesting flexibility. }\\[3mm]
\end{scriptsize}
\section{Introduction}
The purpose of this article is to study the uniqueness of weak solutions for the 3D incompressible  stationary Navier-Stokes equations in $\Rt$. Let us recall that this system is given by the following equations
\begin{equation}\label{NS_section_Stationnaire}
\begin{cases}
\Delta \vu-(\vu \cdot \vn) \vu-\vn P=0,
\\[5pt]
div(\vu)=0,
\end{cases}
\end{equation} 
where $\vu:  \mathbb{R}^3 \longrightarrow \mathbb{R}^3$ is the velocity of the fluid and $P: \mathbb{R}^3 \longrightarrow \mathbb{R}$ denotes its pressure. Although it is an easy exercise to obtain solutions $(\vu, P)$ for this equation in the space $\big(\dot{H}^1(\mathbb{R}^3), \dot{H}^{\frac{1}{2}}(\Rt)\big)$ (see for instance \cite[Theorem 16.2]{lemarie2018navier}), the uniqueness of such solutions is still a challenging open problem. More precisely, we have the following problem (which was initially mentioned in the book \cite[Remark X.9.4 and Theorem X.9.5]{galdi2011introduction} and also stated in the article \cite{Ser2016}):\\ 

\noindent\emph{Show that any solution $\vu$ of the problem (\ref{NS_section_Stationnaire}) which satisfies the conditions
\begin{equation}\label{Conjecture}
\vu \in \dot{H}^1(\R^3) \qquad\mbox{and}\qquad \vu(x)\to 0 \mbox{ as } |x|\to +\infty,
\end{equation}
is identically equal to zero.}\\ 

Remark that, in the setting stated in the problem (\ref{Conjecture}) above, by the classical Sobolev embeddings we have the space inclusion $\dot{H}^1(\Rt)\subset L^6(\Rt)$ from which we can deduce some specific decay at infinity for a solution $\vu$, however this $L^6$-decay seems not enough to conclude that a solution $\vu\in \dot{H}^1(\Rt)$ of the equation (\ref{NS_section_Stationnaire}) is null.  Nevertheless, if we assume an additional hypotheses, for example $\vu\in E(\Rt)$ where $E$ is a ``nice'' functional space (\emph{i.e.} with a stronger decay at infinity than $L^6$), then statements of the following form have been shown:
\begin{center}
\emph{If $\vu\in  \dot{H}^1(\Rt) \cap E(\Rt)$ is a solution of the equation (\ref{NS_section_Stationnaire}) in $\mathbb R^3$, then we have $\vu\equiv 0$,}
\end{center}
and this type of result is known in the literature as a \emph{Liouville-type theorem} for the Navier-Stokes equations (see \emph{e.g.}  \cite{chae2014liouville}, \cite{chae2013liouville}, \cite{ChaeWolf} and \cite{koch2009liouville}). For example, in \cite[Theorem X.9.5]{galdi2011introduction} the case $E=L^{\frac{9}{2}}(\Rt)$ was studied while the space $E=BMO^{-1}(\Rt)$ was considered in \cite{seregin2016liouville}. Other funcional spaces can also be taken into account, see for example the articles \cite{chamorro2021some} and \cite{Kozono}.\\

Let us remark that we can consider a larger functional framework than $\dot{H}^1(\Rt)$ and we can work with weak solutions such that $(\vu,P)\in \big(L^{2}_{loc}(\Rt), \mathcal{D}^{'}(\Rt)\big)$. It is clear that the trivial solution $\vu=0$ satisfies (\ref{NS_section_Stationnaire}), but in this $L^{2}_{loc}$ setting this solution is not unique: indeed, if we define the function $\psi:\Rt \longrightarrow \mathbb{R}$ by $\psi(x_1,x_2,x_3)=\frac{x^{2}_{1}}{2}+\frac{x^{2}_{2}}{2}-x^{2}_{3}$ and if we set the functions $\vu$ and $P$ by  the identities $\vu(x_1,x_2,x_3)=\vec{\nabla}\psi(x_1,x_2,x_3)=(x_1, x_2, -2x_3)$ and $P(x_1,x_2,x_3)=-\frac{1}{2}\vert \vu(x_1,x_2,x_3)\vert^2$, then we have $\vu \in L^{2}_{loc}(\Rt)$ (since $\vert \vu(x)\vert \approx \vert x \vert$) and using the usual rules of vector calculus we have that the couple $(\vu,P)$ given by the expressions above satisfies (\ref{NS_section_Stationnaire}).\\ 

Note now that if we assume an additional information, for example if $E$ is a Lebesgue space $L^p(\Rt)$ with $3\leq p\leq \frac{9}{2}$, then is it possible to prove in the setting $L^{2}_{loc}\cap E$ that the trivial solution is unique (see \cite{chamorro2021some}). However the case when $E=L^p(\Rt)$ with $\frac{9}{2}<p\leq 6$ remains widely open and we can state the following problem:
\begin{equation}\label{Open_Problem}
\begin{split}
&\mbox{\emph{Prove that if $\vu\in L^2_{loc}(\Rt)$ is a weak solution of the stationary Navier-Stokes equations (\ref{NS_section_Stationnaire})}}\\
&\mbox{\emph{such that $\vu\in L^p(\Rt)$ with $\frac{9}{2}<p$, then we have $\vu \equiv 0$.}}\end{split}
\end{equation}
The main feature of this work is thus to explore a uniqueness result for equations (\ref{NS_section_Stationnaire}) using some additional information for the velocity field $\vu$ stated in terms of a Lebesgue spaces of variable exponent $L^{p(\cdot)}(\Rt)$ -instead of the classical Lebesgue spaces $L^p(\Rt)$- and to give some insights to the previous problem (\ref{Open_Problem}).\\

These functional spaces $L^{p(\cdot)}$ are quite different from the usual Lebesgue spaces $L^p$. Indeed, to define the space  $L^{p(\cdot)}(\R^3) $, we will proceed as follows: given a function $p(\cdot):\mathbb{R}^3\longrightarrow [1,+\infty]$, we say that $p(\cdot )  \in \mathcal{P}(\R^3)$ if $p(\cdot )$ is a measurable function.  Then, for a measurable function $\vf:\mathbb{R}^3\longrightarrow \mathbb{R}^3$, we define the \emph{modular function} $\varrho_{p(\cdot)}$ associated to $p(\cdot)$ by the expression
\begin{equation}\label{Def_Modular_Intro}
\varrho_{p(\cdot)}(\vf)=\int_{\mathbb{R}^3}|\vf(x)|^{p(x)}dx.
\end{equation}
Next, we consider the \emph{Luxemburg norm} $\|\cdot\|_{L^{p(\cdot)}}$ associated to the modular function $\varrho_{p(\cdot)}$ (see the books \cite{Cruz_Libro} and \cite{Diening_Libro})
\begin{equation}\label{Def_LuxNormLebesgue}
\|\vf\|_{L^{p(\cdot)}}=\inf\{\lambda > 0: \, \varrho_{p(\cdot)}(\vf/\lambda)\leq1\},
\end{equation}
and we define the Lebesgue spaces of variable exponent $L^{p(\cdot)}(\mathbb{R}^3)$ as the set of measurable functions such that the quantity $\|\cdot\|_{L^{p(\cdot)}}$ given above is finite. Of course there are no simple relationships between the spaces $L^{p(\cdot)}(\Rt)$ and $L^{p}(\Rt)$ (see the previously cited books or Section \ref{Secc_Notaciones_Presentaciones} below for more details).\\ 

We will see in this article how to exploit the flexibility of the variable exponent $p(\cdot)$ in order to deduce a uniqueness result of Liouville type. This flexibility will allow us in Theorem \ref{Theoreme_1} below to consider the integrability interval $]\frac{9}{2}, 6[$ (and beyond) in some regions of the space $\mathbb{R}^3$. The interplay of the variable exponent $p(\cdot)$ with some particular regions of the space $\Rt$ will be studied in more detail in Theorem \ref{Theoreme_2} and Theorem \ref{Theoreme_3}.\\ 

Let us remark that, to the best of our knowledge, the use of Lebesgue spaces of variable exponent $L^{p(\cdot)}$ to establish Liouville type theorems for the stationary Navier-Stokes system (\ref{NS_section_Stationnaire}) seems to be new and we hope that the results presented in this article would shed some light to the problems stated in (\ref{Conjecture}) and (\ref{Open_Problem}).\\

The outline of the article is the following. In Section \ref{Secc_Notaciones_Presentaciones} we first present a small review of the main properties of the spaces $L^{p(\cdot)}$ and then we state our main results. Section \ref{Secc_Liouville} is devoted to the proof of the theorems.
\section{Preliminaries and presentation of the results}\label{Secc_Notaciones_Presentaciones}
For $n\geq 1$, let us first consider a function $p:\mathbb{R}^n\longrightarrow [1,+\infty[$, we will say that $p\in \mathcal{P}(\mathbb{R}^n)$ if $p(\cdot)$ is a measurable function and we define $p^-=\underset{x\in \mathbb{R}^n}{\mbox{inf ess}} \; \{p(x)\}$ and $p^+=\underset{x\in \mathbb{R}^n}{\mbox{sup ess}} \; \{p(x)\}$. In order to distinguish between variable and constant exponents, we will always denote exponent functions by $p(\cdot)$, moreover, for the sake of simplicity and to avoid technicalities 
(see \cite[Chapter 3]{Diening_Libro}),
we will always assume here that we have
$$1<p^- \leq p^+<+\infty,$$
unless specifically stated otherwise.\\ 

With these exponents we can define the spaces $L^{p(\cdot)}(\mathbb{R}^n)$ as the set of measurable functions such that the Luxemburg norm in (\ref{Def_LuxNormLebesgue}) is finite. These functional spaces $L^{p(\cdot)}(\mathbb{R}^n)$ possess  the structural properties of normed spaces and they are moreover Banach spaces, see \cite[Theorem 3.2.7]{Diening_Libro} for more details. \\

Similarly to the classical case, given $p(\cdot)\in \mathcal{P}(\R^3)$  we can define the variable conjugate exponent $q(\cdot)\in \mathcal{P}(\R^3)$ by the pointwise relationship
\begin{equation*}
q(x)= \frac{p(x)}{p(x)-1}. 
\end{equation*} 
Thus, from the previous formula and from the definitions of $p^+$, $p^-$, the following relations hold
\begin{equation}\label{sep_11}
q^- = \frac{p^+}{p^+ -1} \leq \frac{p^-}{p^- -1}= q^+.
\end{equation}
In this setting, the H\"older inequalities have the following version: let $p(\cdot),\,q(\cdot),\,r(\cdot)\in \mathcal{P}(\mathbb{R}^n)$ be functions such that we have the pointwise relationship
$\frac{1}{p(x)}=\frac{1}{q(x)}+\frac{1}{r(x)}$, $x\in \mathbb{R}^n$. Then there exists a constant $C>0$ such that for all $f\in L^{q(\cdot)}(\mathbb{R}^n) $ and $g \in L^{r(\cdot)}(\mathbb{R}^n)$, the pointwise product $fg$ belongs to the space $L^{p(\cdot)}(\mathbb{R}^n)$ and we have the estimate
$$\|f g\|_{L^{p(\cdot)}} \leq C\|f\|_{L^{q(\cdot)}}\|g\|_{L^{r(\cdot)}},$$
see \cite[Theorem 2.26]{Cruz_Libro} or \cite[Lemma 3.2.20]{Diening_Libro} for a proof. This estimate can be easily generalized to vector fields $\vf, \vg :\mathbb{R}^n\longrightarrow \mathbb{R}^n$ and to the product $\vf\cdot \vg$.\\

Now, let us consider a measurable domain $\Omega$ such that $ \Omega \subset \Rt$ and let $p(\cdot) \in \mathcal{P}(\Rt)$. We will denote by $p_\Omega(\cdot)$ the variable exponent restricted to the set $\Omega$, \emph{i.e.} $p_\Omega(\cdot)=  {p} (\cdot)_ { |_{\Omega} }$ and we write
$$p_{\Omega}^-=\underset{x\in \Omega}{\mbox{inf ess}} \; \{p(x)\}\qquad \mbox{and}\qquad p_{\Omega}^+=\underset{x\in \Omega}{\mbox{sup ess}} \; \{p(x)\}.$$
Moreover, by the definition of the Luxemburg norm for $f\in L^{p(\cdot)}(\Rt) $ we have the formula
\begin{equation}\label{sept 15 _1}
\|f\|_{L^{p_\Omega(\cdot)}(\Omega)}=\|f \mathds{1}_{\Omega}\|_{L^{p(\cdot) }(\Rt)},
\end{equation}
see \cite[Section 2.3, page 21]{Cruz_Libro}. The following results will be useful in the sequel:
\begin{Lemme}\label{Proposition_Lp_plus_minus}
Consider a measurable set $\Omega\subset \R^3$ and $p(\cdot)\in \mathcal{P}(\Rt)$ a variable exponent, assume that we have $|\Omega|<+\infty$.  Then 
$$\|1\|_{L^{p_{\Omega}(\cdot) }(\Omega)}\leq 2\max\{|\Omega|^{\frac 1{ p^-}},|\Omega|^{\frac 1{p^+}}\}.$$
\end{Lemme}
The proof of this result can be consulted in \cite[Lemma 3.2.12]{Diening_Libro}. Here is another useful property:
\begin{Lemme}\label{Lemme_Linfty_Lp_variable}
Let $\Omega\subset \R^3$ and $p(\cdot)\in \mathcal{P}(\mathbb{R}^3)$ a variable exponent. Then, we have the space inclusion $L^\infty (\Omega)\subset L^{p_{\Omega}(\cdot)} (\Omega) $, if and only if $1\in  L^{p_{\Omega}(\cdot)} (\Omega) $ and we have the estimate
$$\|f\|_{L^{p_{\Omega}(\cdot)} (\Omega)}\leq \|f\|_{L^\infty(\Omega)}\|1\|_{L^{p_{\Omega}(\cdot)} (\Omega)}.$$
In particular, the embedding holds if $|\Omega|<+\infty$.
\end{Lemme}
The proof of this result can be founded in \cite[Proposition 2.43]{Cruz_Libro}. For more details on the Lebesgue spaces of variable exponent, on their inner structure as well as many other properties, see the books \cite{Cruz_Libro} and \cite{Diening_Libro}. \\

With these preliminaries at our disposal we can state our main results. 
\begin{Theoreme}\label{Theoreme_1} 
Let $(\vu, P)\in (L^2_{loc}(\mathbb{R}^3, \mathcal{D}'(\Rt))$ be a weak solution of the stationary Navier-Stokes equation \eqref{NS_section_Stationnaire}. Consider the following cylinder 
$$\mathscr{C}=\{(x_1,x_2,x_3)\in \R^3: x_2^2  +  x_3^2 \leq 1, \ x_1\in \R \}.$$
Moreover, consider the variable exponent $p(\cdot)\in \mathcal{P}(\R^3)$ defined by
\begin{equation}\label{Def_VariableExpTheo1}
p(x)=
\begin{cases}
3<  p_{(\Rt\setminus \mathscr{C})}^- \leq  p_{(\Rt\setminus \mathscr{C})}(x)\leq  p_{(\Rt\setminus \mathscr{C})}^+< \frac{9}{2},\\[5mm]
\frac{9}{2}< p_{\mathscr{C}}^- \leq  p_{ \mathscr{C}}(x)\leq  p_{ \mathscr{C}}^+<+ \infty.
\end{cases}
\end{equation}
If we assume the additional hypotheses $\vu\in L^{p(\cdot)}(\mathbb{R}^3)$ and $P\in  L^{\frac{p(\cdot)}{2}}(\mathbb{R}^3)$, then we have $\vu=0$.
\end{Theoreme}
Some remarks are in order here. First note that, besides the additional condition $\vu\in L^{p(\cdot)}(\mathbb{R}^3)$ over the velocity field, we also demand the hypothesis $P\in  L^{\frac{p(\cdot)}{2}}(\mathbb{R}^3)$ for the pressure. Let us explain this constraint: recall that, using the divergence-free property of $\vu$, we have the classical relationship $-\Delta P=div(div(\vu\otimes \vu))$ from which we deduce the expression $\displaystyle{P=\sum_{i,j=1}^3\mathcal{R}_i\mathcal{R}_i(u_iu_j)}$ where $\mathcal{R}_i$ are the usual Riesz transforms. We can thus obtain quite easily some information over the pressure $P$ from the information available over $\vu$, \emph{as long as} the Riesz transforms are bounded in the functional framework considered. However, to the best of our knowledge, in the case of the Lebesgue spaces of variable exponents, the Riesz transforms \emph{are not bounded} in $L^{p(\cdot)}(\Rt)$ unless the exponent $p(\cdot)$ satisfies the lower and upper bounds $1<p^-\leq p(\cdot)\leq p^+<+\infty$ and one additional condition (see \cite[Lema 12.4.3]{Diening_Libro}) which can be easily stated in terms of the log-H\"older regularity property\footnote{For the sake of completeness, we point out that the boundedness of the Riesz transforms in the spaces $L^{p(\cdot)}(\Rt)$ can be expressed by a more general and more technical condition: we have $\|\mathcal{R}_i(f)\|_{L^{p(\cdot)}(\Rt)}\leq C\|f\|_{L^{p(\cdot)}(\Rt)}$ if we have $1<p^-\leq p^+<+\infty$ and $p(\cdot)\in \mathcal{A}$, where the rather technical definition of the class of variable exponents $\mathcal{A}$ is given in \cite[Definition 4.4.6]{Diening_Libro}. Of course we have that if $p(\cdot)\in \mathcal{P}^{log}(\Rt)$ then we have $p(\cdot)\in \mathcal{A}$. See the book \cite{Diening_Libro} for more details.}: we will say that $p(\cdot)$ is log-H\"older continuous (denoted by $p(\cdot)\in \mathcal{P}^{log}(\Rt)$) if
\begin{eqnarray}
\left|\frac{1}{p(x)}-\frac{1}{p(y)}\right|&\leq &\frac{C}{\log(e+1/|x-y|)}\qquad \mbox{for all } x,y\in \mathbb{R}^3,\quad \mbox{and}\notag\\
 \left|\frac{1}{p(x)}-\frac{1}{p_\infty}\right|&\leq&\frac{C}{\log(e+|x|)}\qquad \mbox{for all } x\in \mathbb{R}^3,\quad \mbox{where} \quad \frac{1}{p_\infty}=\underset{|x|\to +\infty}{\lim}\frac{1}{p(x)}.\label{LogHolder}
\end{eqnarray}
See Definition 4.1.1 and Definition 4.1.4 of the book \cite{Diening_Libro} for more details about the log-H\"older regularity property. Remark that the condition (\ref{LogHolder}) above imposes a specific behavior of the variable exponents at infinity which is not satisfied by the function $p(\cdot)$ considered in (\ref{Def_VariableExpTheo1}). Thus, since we can not deduce as easily as for the classical Lebesgue spaces some interesting control over the pressure $P$ from the information of $\vu$, we therefore ask here the constraint $P\in  L^{\frac{p(\cdot)}{2}}(\mathbb{R}^3)$ which does not interfere with the general purposes of our theorems.\\

Next, we remark that over the set $\Rt\setminus \mathscr{C}$ (\emph{i.e.} outside the cylinder $\mathscr{C}$) we impose to the variable exponent $p(\cdot)$ to be in the open interval $]3, \frac{9}{2}[$ and these are the known restrictions from which we can deduce the uniqueness of the trivial solution (see \cite{chamorro2021some}). The novelty comes then from the conditions \emph{inside} the cylinder $\mathscr{C}$ where we can consider integrability values beyond this range, in fact the constraint (\ref{Def_VariableExpTheo1}) was precisely made to study case $p(\cdot)>\frac{9}{2}$ and this result gives some information related to the problem (\ref{Open_Problem}). Note also that the Lebesgue measure of the cylinder is infinite (\emph{i.e.} $|\mathscr{C}|=+\infty$) and the previous result can be extended to many different subsets of $\Rt$ as long as they satisfy some simple restrictions (for example a finite union of such cylinders): this shows that we can obtain the uniqueness of the trivial solution with a slow decay at infinity in some non-negligible regions of the space $\Rt$. \\

Finally, we note that in the definition of the variable exponent $p(\cdot)$ given (\ref{Def_VariableExpTheo1}) we are considering (for the sake of simplicity) the lower bound $3<p^-$ and open intervals in the different constraints. Some of these conditions may probably be relaxed and we do not claim any kind of optimality in our results.\\

In the following theorem, we study in more detail the behavior of the variable exponent $p(\cdot)$ in the case $\frac{9}{2}<p(\cdot)$ over larger subsets of $\Rt$. 

\begin{Theoreme}\label{Theoreme_2} 
Let $(\vu, P)\in (L^2_{loc}(\mathbb{R}^3, \mathcal{D}'(\Rt))$ be a weak solution of the stationary Navier-Stokes equation \eqref{NS_section_Stationnaire}. Consider now $\mathscr{S}$ the set defined by
$$\mathscr{S}=\{(x_1, x_2, x_3)\in \Rt: x_2^2+x_3^2\leq x_1^\gamma, \; x_1>0\},$$
where $0<\gamma<1$. We define the variable exponent $\mathfrak{p}(\cdot)$ by the following conditions:
\begin{equation}\label{Def_VariableExpTheo2}
\mathfrak{p}(x)=
\begin{cases}
3<  \mathfrak{p}_{(\Rt\setminus \mathscr{S})}^- \leq  \mathfrak{p}_{(\Rt\setminus \mathscr{S})}(x)\leq  \mathfrak{p}_{(\Rt\setminus \mathscr{S})}^+< \frac{9}{2},\\[5mm]
\frac{9}{2}< \mathfrak{p}_{\mathscr{S}}^- \leq  \mathfrak{p}_{ \mathscr{S}}(x)\leq  \mathfrak{p}_{ \mathscr{S}}^+< \frac{6\gamma+3}{2\gamma}.
\end{cases}
\end{equation}
If we have the restriction $\vu\in L^{\mathfrak{p}(\cdot)}(\mathbb{R}^3)$ and $P\in L^{\frac{\mathfrak{p}(\cdot)}{2}}(\mathbb{R}^3)$, then we obtain that $\vu=0$.
\end{Theoreme}
As we can observe here, the larger the set $\mathscr{S}$ is (reflected by $\gamma\to 1^-$), the more the upper bound for the variable exponent $\mathfrak{p}_{ \mathscr{S}}^+$ tends towards known values (\emph{i.e.} we have $\frac{6\gamma+3}{2\gamma}\to \frac{9}{2}^+$). Note also that, if $\gamma\to 0^+$ the shape of the set $\mathscr{S}$ tends to the cylinder $\mathscr{C}$ considered in Theorem \ref{Theoreme_1} above and then we have that the upper bound satisfy $\frac{6\gamma+3}{2\gamma}\to +\infty$, which is the second condition stated in the expression (\ref{Def_VariableExpTheo1}). This theorem shows how the variable exponent $\mathfrak{p}(\cdot)$ may vary over a set $\mathscr{S}$ whose Lebesgue measure is also variable: indeed, if the measure of this set is ``reasonable'' then we can consider a mild decay at infinity and we can still obtain the uniqueness of the trivial solution. However, if the measure of this set is ``too big'' then we shall recover the known upper bound $\frac{9}{2}$ which allows us to solve the uniqueness problem considered here. Note that when $\gamma=1$ this result gives an example of such \emph{big} sets. \\

It is worth to remark now that, although integrability values outside the interval $]3, \frac{9}{2}[$ can be considered in the two previous results (namely $\frac{9}{2}<p(\cdot), \mathfrak{p}(\cdot)$ over suitable subsets of $\Rt$, $\mathscr{C}$ or $\mathscr{S}$), we always have the upper bounds $p^+, \mathfrak{p}^+<+\infty$. In the next theorem we study the case when we can have $p(\cdot)=+\infty$ over a particular subset of $\Rt$.
\begin{Theoreme}\label{Theoreme_3} 
Let $(\vu, P)\in (L^2_{loc}(\mathbb{R}^3, \mathcal{D}'(\Rt))$ be a weak solution of the stationary Navier-Stokes equation \eqref{NS_section_Stationnaire}. We define the set $\mathscr{N}$ by the condition
$$\mathscr{N}=\{(x_1, x_2, x_3)\in \Rt: x_2^2+x_3^2\leq x_1^{-\frac{\sigma}{2}},\;\; x_1>0\},$$
with $0<\sigma<1$. Consider now the variable exponent ${\bf p}(\cdot)$ given by the following conditions:
\begin{equation}\label{Def_VariableExpTheo3}
{\bf p}(x)=
\begin{cases}
3<  {\bf p}_{(\Rt\setminus \mathscr{N})}^- \leq  {\bf p}_{(\Rt\setminus \mathscr{N})}(x)\leq  {\bf p}_{(\Rt\setminus \mathscr{N})}^+< \frac{9}{2},\\[5mm]
{\bf p}_{ \mathscr{N}}(x)=+\infty.
\end{cases}
\end{equation}
If we have the restriction $\vu\in L^{{\bf p}(\cdot)}(\mathbb{R}^3)$ and $P\in L^{\frac{{\bf p}(\cdot)}{2}}(\mathbb{R}^3)$, then we obtain that $\vu=0$.
\end{Theoreme}
First note that the Lebesgue measure of the set $\mathscr{N}$ is infinite, so we can consider here vector fields $\vu$ satisfying the equation (\ref{NS_section_Stationnaire}) which are $L^2_{loc}(\Rt)$, have a suitable decay at infinity in a large portion of the space $\Rt$ but that can be constant at infinity over a non-negligible subset of $\Rt$. Thus considering all these restrictions, we can deduce the uniqueness of the trivial solution and this result may suggest that some decay at infinity is not absolutely necessary to solve the problem (\ref{Conjecture}). Of course, and as pointed out before, we do not claim any optimality in our results and many others subsets of the space $\Rt$ can be studied with similar conclusions. 
\section{Proof of the Theorems}\label{Secc_Liouville}
The proof of the three previous theorems is relatively similar: we start by localizing the information with a smooth cut-off function supported in a ball $B(0,R)$ and then we analyze the behavior of the localized information as $R\to +\infty$: it is in this step that we will exploit the different hypotheses stated in Theorems \ref{Theoreme_1}, \ref{Theoreme_2} and \ref{Theoreme_3} to deduce the uniqueness of the trivial solution.
\subsection{General Framework}\label{SubsecGeneral}
Let us start with a divergence-free vector field $\vu\in L^2_{loc}(\Rt)$ that satisfies in the weak sense the stationary Navier-Stokes equations (\ref{NS_section_Stationnaire}). Since by \cite[Theorem 2.51]{Cruz_Libro} we have the space inclusion 
$$L^{p(\cdot)}(\Rt)\subset L^{p^-}(\Rt)+L^{p^+}(\Rt),$$
and since by the definition of the variable exponent $p(\cdot)$ (or $\mathfrak{p}(\cdot)$ or ${\bf p}(\cdot)$) used in our theorems we always have the lower and upper bounds $3<p^-\leq p^+\leq +\infty$ (see (\ref{Def_VariableExpTheo1}) for $p(\cdot)$, (\ref{Def_VariableExpTheo2}) for $\mathfrak{p}(\cdot)$ and (\ref{Def_VariableExpTheo3}) for ${\bf p}(\cdot)$), then we deduce the following embeddings 
$$L^{p(\cdot)}(\Rt)\subset L^{p^-}(\Rt)+L^{p^+}(\Rt)\subset L^{p^-}_{loc}(\Rt)+L^{p^+}_{loc}(\Rt)\subset L^3_{loc}(\Rt).$$
Now, as we have $\vu\in L^3_{loc}(\Rt)$, then following \cite[Theorem X.1.1]{galdi2011introduction} we obtain that the velocity field $\vu$ and the pressure $P$ are regular functions.\\

Now, let $\theta\in {\mathcal{C}}^\infty_0 (\mathbb{R}^3) $ be a cut-off function such that $0<\theta<1$, $\theta(x)=1$ if $|x|<1$, $\theta(x)=0$ if $|x|>2$. Given $R>1$, we define the function $\theta_R$ by $\theta_R(x)=\theta(\frac{x}{R})$: thus, $\theta_R(x)=1$ if $|x|< \frac R 2$ and $\theta_R(x)=0$ if $|x|\geq R$.
By testing the Navier-Stokes equations \eqref{NS_section_Stationnaire} with $\theta_R \vu$ and using the fact that $supp(\theta_R \vu)=B_R=B(0,R)$ we obtain 
\begin{equation}\label{eq. 3}
\int_{B_R}-\Delta \vu \cdot\left(\theta_R \vu\right)+(\vu \cdot \vn) \vu \cdot\left(\theta_R \vu\right)+\vn P \cdot\left(\theta_R \vu\right) d x=0.
\end{equation}
Note that since $\vu$ and $P$ are smooth enough, all the terms in the identity above are well-defined. Then, by using the fact that $div(\vu)=0$ and integrating by parts, we obtain that the terms in the left-hand side of \eqref{eq. 3} can be rewritten respectively as 
\begin{eqnarray*}
\int_{B_R}-\Delta \vu \cdot\left(\theta_R \vu\right) dx&=&-\int_{B_R} \Delta \theta_R\left(\frac{|\vu|^2}{2}\right) d x+\int_{B_R} \theta_R|\vn \otimes u|^2 dx,\\
\int_{B_R}(\vu \cdot \vn) \vu \cdot\left(\theta_R \vu\right) dx&=&\sum_{i, j=1}^3 \int_{\Rt} \theta_R \partial_j\left(u_j \frac{u_i^2}{2}\right) d x-\int_{B_R}\vn\theta_R \cdot\left(\frac{|\vu|^2}{2} \vu\right) d x,\\
\int_{B_R} \vn P \cdot\left(\theta_R \vu\right) d x &=&-\int_{B_R} \vn \theta_R \cdot(P \vu)dx.
\end{eqnarray*}
Considering these identities, we can rewrite \eqref{eq. 3} in the following manner
\begin{equation*}
\int_{B_R} \theta_R|\vn \otimes \vu|^2 dx=\int_{B_R} \Delta \theta_R \frac{|\vu|^2}{2} dx+
\int_{B_R} \vn \theta_R \cdot\left(\left(\frac{|\vu|^2}{2}+P\right) \vu\right)dx.
\end{equation*}
Since we have $\theta_R(x)\equiv 1$ over the set $|x|< \frac R 2$, we can write
\begin{equation}\label{ineq_Base}
\int_{B_{\frac{R}{2}}}|\vn \otimes \vu|^2 dx \leq \underbrace{\int_{B_R} \Delta \theta_R \frac{|\vu|^2}{2} dx}_{\alpha(R)}+\underbrace{\int_{B_R} \vn \theta_R \cdot\left(\left(\frac{|\vu|^2}{2}+P\right) \vu\right)dx}_{\beta(R)}. 
\end{equation}
To conclude that we have $\vu=0$ in the context of Theorems \ref{Theoreme_1}, \ref{Theoreme_2} and \ref{Theoreme_3}, we aim to prove in the next subsections that we have the following limits
\begin{equation}\label{LimitsAlphaBeta}
\displaystyle \lim_{R\to +\infty}  \alpha(R) =\lim_{R\to +\infty}  \beta(R)=0.
\end{equation}
Indeed, if we establish these limits, we can conclude from the estimate (\ref{ineq_Base}) above that 
\begin{equation}\label{LimiteSobolev}
\lim_{R\to +\infty}\int_{B_{\frac{R}{2}}}|\vn \otimes \vu|^2 dx=\|\vu\|_{\dot{H}^1}=0,
\end{equation}
from which we deduce, by the Sobolev embeddings, that $\|\vu\|_{L^6}=0$ and thus that we have $\vu\equiv 0$. 
\subsection{Proof of the Theorem \ref{Theoreme_1}} 
As said previously, we only need now to study the terms $\alpha(R)$ and $\beta(R)$ defined in the expression (\ref{ineq_Base}) and to show that these quantities tend to $0$ as $R\to+\infty$.
\begin{itemize}
\item[{\bf 1)}]{\bf Control for $\alpha(R)$}. For studying the term $\alpha(R)$ in \eqref{ineq_Base}, the H\"older inequality\footnote{Note that by (\ref{Def_VariableExpTheo1}) we always have $1\leq \frac{p(\cdot)}{2}\leq +\infty$ and thus we can apply the H\"older inequalities in the setting of Lebesgue spaces of variable exponents.}  with $\frac{2}{p(\cdot)}+\frac{1}{q(\cdot)}=1$ (see \cite[Theorem 2.26]{Cruz_Libro}) yields the following estimate 
$$\alpha(R)=\int_{B_R} \Delta \theta_R \frac{|\vu|^2}{2} dx \leq \int_{\Rt} \Delta \theta_R \frac{|\vu|^2}{2} dx\leq C\| \Delta \theta_R \|_{L^{q(\cdot)}(\Rt)}\| |\vu|^2 \|_{L^{\frac{p(\cdot)}{2}}(\Rt)},$$
since we have the identity $\| |\vu|^2 \|_{L^{\frac{p(\cdot)}{2}}(\Rt)}=\|\vu \|_{L^{p(\cdot)}(\Rt)}^2$ (see \cite[Proposition 2.18]{Cruz_Libro}) we can write
\begin{equation}\label{18aveq32_T1}
\alpha(R)\leq C\|\Delta \theta_R \|_{L^{q(\cdot)}(\Rt)}\|\vu\|_{L^{p(\cdot)}(\Rt)}^2.
\end{equation}
Note that, by hypothesis we have $\|\vu\|_{L^{p(\cdot)}(\Rt)}<+\infty$, so we only need to estimate the quantity $\|\Delta \theta_R \|_{L^{q(\cdot)}(\Rt)}$. Since the variable exponent $p(\cdot)$ satisfies the relationships stated in (\ref{Def_VariableExpTheo1}), from the pointwise H\"older relationship $\frac{2}{p(\cdot)}+\frac{1}{q(\cdot)}=1$, we deduce that the variable exponent $q(\cdot)$ satisfies the equation $q(\cdot)=\frac{p(\cdot)}{p(\cdot)-2}$, \emph{i.e.}: 
\begin{equation}\label{Definition_qConjugueTheo1}
q(x)=
\begin{cases}
\frac{9}{5}<  q_{(\Rt\setminus \mathscr{C})}^- \leq  q_{(\Rt\setminus \mathscr{C})}(x)\leq  q_{(\Rt\setminus \mathscr{C})}^+< 3,\\[5mm]
1< q_{\mathscr{C}}^- \leq  q_{ \mathscr{C}}(x)\leq  q_{ \mathscr{C}}^+<\frac{9}{5},
\end{cases}
\end{equation}
where $\mathscr{C}=\{(x_1,x_2,x_3)\in \R^3: x_2^2  +  x_3^2 \leq 1, \ x_1\in \R \}$.\\

Now, to estimate the first term in the right-hand side of the expression (\ref{18aveq32_T1}) above, 
we recall that from the definition of the function $\theta_R$ we have $supp \left( \Delta \theta_R\right) \subset \mathcal{C}(\frac{R}{2},R)=\{x\in \Rt: \frac{R}{2}\leq |x|\leq R\}$ and we have
\begin{equation}\label{ConditionSupportTheta}
\|\Delta\theta_R\|_{L^{q(\cdot)}(\Rt)}=\|\Delta\theta_R\|_{L^{q(\cdot)}(\mathcal{C}(\frac{R}{2},R))}.
\end{equation}
To continue, we denote by $\mathscr{C}_1$ and $\mathscr{C}_2$ the sets defined by 
\begin{equation}\label{Def_SubConjuntos}
\mathscr{C}_1=\mathcal{C}\left(\tfrac{R}{2}, R\right)\cap \mathscr{C}\qquad \mbox{and}\qquad \mathscr{C}_2=\mathcal{C}\left(\tfrac{R}{2}, R\right) \setminus   \mathscr{C},
\end{equation}
respectively. Then, from the property \eqref{sept 15 _1}  we have
\begin{eqnarray}
\| \Delta \theta_R \|_{L^{q(\cdot)}(\mathcal{C}(\frac{R}{2},R))}&=&\| \Delta \theta_R (\mathds{1}_{\mathscr{C}_1}+\mathds{1}_{\mathscr{C}_2}) \|_{L^{q(\cdot)}(\mathcal{C}(\frac{R}{2},R))}\notag\\
&\leq&\|\Delta \theta_R \|_{L^{q_{\mathscr{C}_1}(\cdot)}(\mathscr{C}_1)}+\| \Delta \theta_R \|_{L^{q_{\mathscr{C}_2}(\cdot)}(\mathscr{C}_2)}.\label{sep_1T1}
\end{eqnarray}
We will study these two terms separately.
\begin{itemize}
\item[$\bullet$] In order to control the quantity $\|\Delta \theta_R \|_{L^{q_{\mathscr{C}_1}(\cdot)}(\mathscr{C}_1)}$ above, by the Lemma \ref{Lemme_Linfty_Lp_variable} we can write
\begin{equation}\label{EstimationNormeLinfiniTheta_R1}
\|\Delta \theta_R \|_{L^{q_{\mathscr{C}_1}(\cdot)}(\mathscr{C}_1)}\leq C\|\Delta \theta_R  \|_{L^{\infty}  (\mathscr{C}_1)}\|1\|_{L^{q_{\mathscr{C}_1}(\cdot)}(\mathscr{C}_1)}.
\end{equation}
Then, by the definition of the function $\theta_R$, since $\mathscr{C}_1=\mathcal{C}(\frac{R}{2}, R)\cap \mathscr{C}\subset \mathcal{C}(\frac{R}{2}, R)$, we have $\|\Delta \theta_R  \|_{L^{\infty}  (\mathscr{C}_1)}\leq \|\Delta \theta_R  \|_{L^{\infty}(\mathcal{C}(\frac{R}{2}, R))}\leq CR^{-2}$ and we obtain
$$\|\Delta \theta_R \|_{L^{q_{\mathscr{C}_1}(\cdot)}(\mathscr{C}_1)}\leq C R^{-2}\|1\|_{L^{q_{\mathscr{C}_1}(\cdot)}(\mathscr{C}_1)}.$$
Now, as $\mathscr{C}_1=\mathcal{C}(\frac{R}{2}, R)\cap \mathscr{C}\subset \mathscr{C}$, we have $q^-_{\mathscr{C}}\leq q^-_{\mathscr{C}_1}\leq q^+_{\mathscr{C}_1}\leq q^+_{\mathscr{C}}$ and thus, applying the Lemma \ref{Proposition_Lp_plus_minus} to estimate the quantity $\|1\|_{L^{q_{\mathscr{C}_1}(\cdot)}(\mathscr{C}_1)}$, we can write
\begin{equation}\label{EstimationC01}
\|\Delta \theta_R \|_{L^{q_{\mathscr{C}_1}(\cdot)}(\mathscr{C}_1)}\leq C R^{-2}\max \{|\mathcal{C}(\tfrac{R}{2}, R)\cap \mathscr{C}|^{\frac 1 {q_{\mathscr{C}}^- }}, |\mathcal{C}(\tfrac{R}{2}, R)\cap \mathscr{C}|^{\frac 1 {q_{\mathscr{C}}^+}}\}.
\end{equation}
Noticing that we have the set inclusion 
$$\mathcal{C}(\tfrac{R}{2}, R)\cap \mathscr{C}\subset \mathscr{A}=\{(x_1,x_2,x_3)\in \R^3: x_2^2  +  x_3^2 \leq 1, -R< x_1<R \},$$ 
we easily deduce that $|\mathcal{C}(\tfrac{R}{2}, R)\cap \mathscr{C}|\leq |\mathscr{A}|=CR$, since $\mathscr{A}$ is a cylinder of diameter $1$ and height $2R$. We thus obtain
\begin{equation}\label{EstimationC1}
\|\Delta \theta_R \|_{L^{q_{\mathscr{C}_1}(\cdot)}(\mathscr{C}_1)}\leq C\max \{R^{-2+\frac 1{q_{\mathscr{C}}^-}}, R^{-2+\frac 1 {q_{\mathscr{C}}^+}}\},
\end{equation}
but since, by (\ref{Definition_qConjugueTheo1}), we have $\frac{5}{9}<\frac 1 {q_{\mathscr{C}}^+}\leq \frac 1 {q_{\mathscr{C}}^-}<1$, it comes $-2+\frac 1 {q_{\mathscr{C}}^+}\leq -2+\frac 1 {q_{\mathscr{C}}^-}<-1$ from which we obtain that 
\begin{equation}\label{Limite_Theo111}
\|\Delta \theta_R \|_{L^{q_{\mathscr{C}_1}(\cdot)}(\mathscr{C}_1)}\underset{R\to +\infty}{\longrightarrow}0.
\end{equation}
\item[$\bullet$] We consider now the second term in the right-hand side of (\ref{sep_1T1}). By the same arguments as above we can write
\begin{equation}\label{EstimationNormeLinfiniTheta_R2}
\|\Delta \theta_R \|_{L^{q_{\mathscr{C}_2}(\cdot)}(\mathscr{C}_2)}\leq C\|\Delta \theta_R\|_{L^{\infty}  (\mathscr{C}_2)}\|1\|_{L^{q_{\mathscr{C}_2}(\cdot)}(\mathscr{C}_2)},
\end{equation}
again, since $\mathscr{C}_2=\mathcal{C}(\frac{R}{2}, R)\setminus \mathscr{C}\subset \mathcal{C}(\frac{R}{2}, R)$ we have $\|\Delta \theta_R  \|_{L^{\infty}  (\mathscr{C}_2)}\leq \|\Delta \theta_R  \|_{L^{\infty}(\mathcal{C}(\frac{R}{2}, R))}\leq CR^{-2}$, and we obtain
$$\|\Delta \theta_R \|_{L^{q_{\mathscr{C}_2}(\cdot)}(\mathscr{C}_2)}\leq CR^{-2}\|1\|_{L^{q_{\mathscr{C}_2}(\cdot)}(\mathscr{C}_2)},$$
and applying Lemma \ref{Proposition_Lp_plus_minus} we deduce the estimate
$$\|\Delta\theta_R\|_{L^{q_{\mathscr{C}_2}(\cdot)}(\mathscr{C}_2)}\leq CR^{-2}\max\{|\mathscr{C}_2|^{\frac 1{q_{\mathscr{C}_2}^-}},|\mathscr{C}_2|^{\frac 1{q_{\mathscr{C}_2}^+}}\}.$$
At this point we remark that we have the set inclusions $\mathscr{C}_2=\mathcal{C}(\frac{R}{2}, R)\setminus \mathscr{C}\subset \Rt\setminus \mathscr{C}$ so we obtain, by the first relationship in  (\ref{Definition_qConjugueTheo1}), that $\frac{9}{5}< q_{(\Rt\setminus\mathscr{C})}^-<q_{\mathscr{C}_2}^-\leq q_{\mathscr{C}_2}^+<q_{(\Rt\setminus\mathscr{C})}^+< 3$, \emph{i.e.}
$$\frac{1}{3}< \frac{1}{q_{(\Rt\setminus\mathscr{C})}^+}<\frac{1}{q_{\mathscr{C}_2}^+}\leq \frac{1}{q_{\mathscr{C}_2}^-}<\frac{1}{q_{(\Rt\setminus\mathscr{C})}^-}<  \frac{5}{9}.$$
Moreover, since $|\mathscr{C}_2|=|\mathcal{C}(\frac{R}{2}, R)\setminus \mathscr{C}|\leq |\mathcal{C}(\frac{R}{2}, R)|\leq CR^3$ we have
\begin{eqnarray}
\|\Delta\theta_R\|_{L^{q_{\mathscr{C}_2}(\cdot)}(\mathscr{C}_2)}&\leq &CR^{-2}\max\{R^{\frac3{q_{(\Rt\setminus\mathscr{C})}^-}},R^{\frac 3{q_{(\Rt\setminus\mathscr{C})}^+}}\}\notag\\
&\leq &C\max\{R^{-2+\frac3{q_{(\Rt\setminus\mathscr{C})}^-}},R^{-2+\frac 3{q_{(\Rt\setminus\mathscr{C})}^+}}\}.\label{EstimationC2}
\end{eqnarray}
Observing that we have the bounds $-2+\frac{3}{q_{\Rt\setminus\mathscr{C}}^+}<-2+\frac{3}{q_{\Rt\setminus\mathscr{C}}^-}<  -\frac{1}{3}$, we can deduce that 
\begin{equation}\label{Limite_Theo112}
\|\Delta \theta_R \|_{L^{q_{\mathscr{C}_2}(\cdot)}(\mathscr{C}_2)}\underset{R\to +\infty}{\longrightarrow}0.
\end{equation}
\end{itemize}
With the information (\ref{Limite_Theo111}) and (\ref{Limite_Theo112}), we can come back to the estimate (\ref{sep_1T1}) and we obtain that $\|\Delta \theta_R \|_{L^{q(\cdot)}(\mathcal{C}(\frac{R}{2},R))}=\|\Delta \theta_R \|_{L^{q(\cdot)}(\Rt)}\underset{R\to +\infty}{\longrightarrow}0$. From this control and from the estimate (\ref{18aveq32_T1}), we can directly deduce that $\alpha(R)\underset{R\to +\infty}{\longrightarrow}0$.

\item[{\bf 2)}]{\bf Control for $\beta(R)$}. We recall that from the definition of $\theta_R$ we have $supp( \vn \theta_R)\subset\mathcal{C}(\frac{R}{2},R)$. Thus, by the definition of the term $\beta(R)$ given in (\ref{ineq_Base}) we can write 
\begin{eqnarray}
\left|\beta(R)\right|&=&\left|\int_{B_R} \vn \theta_R \cdot\left(\left(\frac{|\vu|^2}{2}+P\right) \vu\right)dx\right|\notag\\
& \leq &\frac{1}{2} \underbrace{\int_{\mathcal{C}\left(\frac{R}{2}, R\right)}|\vn \theta_R||\vu|^3dx}_{\beta_1(R)}+\underbrace{\int_{\mathcal{C}\left(\frac{R}{2}, R\right)}|\vn \theta_R||P \| \vu| dx}_{\beta_2(R)}.\label{aug_31T1}
\end{eqnarray}
We aim to prove now that we have $\displaystyle \lim_{R\to +\infty}\beta_1(R)= \lim_{R\to +\infty} \beta_2(R)=0$ and these two limits will be studied separately.
\begin{itemize}
\item[$\bullet$] For the term $\beta_1(R)$ we write, by the H\"older inequality\footnote{Recall that by (\ref{Def_VariableExpTheo1}) we always have $1\leq\frac{p(\cdot)}{3}\leq +\infty$ so we can apply the H\"older inequality with parameter $\frac{p(\cdot)}{3}$.} with $\frac{3}{p(\cdot)}+\frac{1}{r(\cdot)}=1$:
$$\beta_1(R)=\int_{\mathcal{C}\left(\frac{R}{2}, R\right)}|\vn \theta_R||\vu|^3dx \leq C\|\vn\theta_R \|_{L^{ r(\cdot)}(\mathcal{C}\left(\frac{R}{2}, R\right))}\| |\vu|^3 \|_{L^{\frac{p(\cdot)}{3}}(\mathcal{C}\left(\frac{R}{2}, R\right))},$$
now, since $\| |\vu|^3 \|_{L^{\frac{p(\cdot)}{3}}(\mathcal{C}\left(\frac{R}{2}, R\right))}\leq \| |\vu|^3 \|_{L^{\frac{p(\cdot)}{3}}(\Rt)}$, we have 
\begin{eqnarray}
\beta_1(R)&\leq &C\|\vn\theta_R \|_{L^{ r(\cdot)}(\mathcal{C}\left(\frac{R}{2}, R\right))}\| |\vu|^3 \|_{L^{\frac{p(\cdot)}{3}}(\Rt)} \notag\\
&\leq &C\|\vn\theta_R \|_{L^{ r(\cdot)}(\mathcal{C}\left(\frac{R}{2}, R\right))}\|\vu \|_{L^{p(\cdot)}(\Rt)}^3,\label{2mayeq3T1}
\end{eqnarray}
where we used the identity $\| |\vu|^3 \|_{L^{\frac{p(\cdot)}{3}}(\Rt)}=\|\vu \|_{L^{p(\cdot)}(\Rt)}^3$. Since, by hypothesis we have $\|\vu \|_{L^{p(\cdot)}(\Rt)}<+\infty$, to prove that $\displaystyle \lim_{R\to +\infty}\beta_1(R)=0$ we only need to study the quantity $\|\vn\theta_R \|_{L^{ r(\cdot)}(\mathcal{C}\left(\frac{R}{2}, R\right))}$. For this, we note that as the variable exponent $p(\cdot)$ satisfies the relationships stated in (\ref{Def_VariableExpTheo1}), from the pointwise H\"older relationship $\frac{3}{p(\cdot)}+\frac{1}{r(\cdot)}=1$, we have $r(\cdot)=\frac{p(\cdot)}{p(\cdot)-3}$, \emph{i.e.}: 
\begin{equation}\label{Definition_rConjugueTheo1}
r(x)=
\begin{cases}
3<  r_{(\Rt\setminus \mathscr{C})}^- \leq  r_{(\Rt\setminus \mathscr{C})}(x)\leq  r_{(\Rt\setminus \mathscr{C})}^+< +\infty,\\[5mm]
1< r_{\mathscr{C}}^- \leq  r_{ \mathscr{C}}(x)\leq  r_{ \mathscr{C}}^+<3,
\end{cases}
\end{equation}
where $\mathscr{C}=\{(x_1,x_2,x_3)\in \R^3: x_2^2  +  x_3^2 = 1, \ x_1\in \R \}$. Using the sets $\mathscr{C}_1$ and $\mathscr{C}_2$ defined in (\ref{Def_SubConjuntos}) and proceeding just as in the estimate (\ref{sep_1T1}) above, we can write
$$\| \vn \theta_R \|_{L^{r(\cdot)}(\mathcal{C}(\frac{R}{2},R))}\leq\|\vn \theta_R \|_{L^{r_{\mathscr{C}_1}(\cdot)}(\mathscr{C}_1)}+\| \vn \theta_R \|_{L^{r_{\mathscr{C}_2}(\cdot)}(\mathscr{C}_2)}.$$
Following the same ideas that leaded us to the estimates (\ref{EstimationC1}) and (\ref{EstimationC2}) (with the difference that $ \|\vn \theta_R  \|_{L^{\infty}}\leq CR^{-1}$) we obtain
\begin{equation}\label{EstimationBeta1C}
\| \vn \theta_R \|_{L^{r(\cdot)}(\mathcal{C}(\frac{R}{2},R))}\leq C\max \{R^{-1+\frac 1{r_{\mathscr{C}}^-}}, R^{-1+\frac 1 {r_{\mathscr{C}}^+}}\}+C\max\{R^{-1+\frac3{r_{(\Rt\setminus\mathscr{C})}^-}},R^{-1+\frac 3{r_{(\Rt\setminus\mathscr{C})}^+}}\},
\end{equation}
now, by the values of the function $r(\cdot)$ given in (\ref{Definition_rConjugueTheo1}) we easily have that 
$$-1+\frac 1 {r_{\mathscr{C}}^+}\leq -1+\frac 1 {r_{\mathscr{C}}^-}<0\qquad \mbox{and}\qquad -1+\frac3{r_{(\Rt\setminus\mathscr{C})}^+}\leq -1+\frac3{r_{(\Rt\setminus\mathscr{C})}^-}<0,$$
from which we deduce that $\underset{R\to +\infty}{\lim}\|\vn \theta_R \|_{L^{r(\cdot)}(\mathcal{C}(\frac{R}{2},R))}=0$, and thus, by the control (\ref{2mayeq3T1}) we have $\displaystyle \lim_{R\to +\infty}\beta_1(R)=0$.

\item[$\bullet$] We continue now with the analysis of the term $\beta_2(R)$ given in the expression (\ref{aug_31T1}). Applying the H\"older inequalities with $\frac{1}{p(\cdot)}+\frac{2}{p(\cdot)}+\frac{1}{r(\cdot)}=1$ (see \cite[Corollary 2.30]{Cruz_Libro}), we  have
\begin{eqnarray*}
\beta_2(R)&=&\int_{\mathcal{C}\left(\frac{R}{2}, R\right)}|\vn \theta_R||P \| \vu| dx\leq C \| \vn\theta_R \|_{L^{r(\cdot)}(\mathcal{C}(\frac{R}{2},R))} \|P \|_{L^{\frac{p(\cdot)}{2}}(\mathcal{C}(\frac{R}{2},R))}\|\vu \|_{L^{p(\cdot)}(\mathcal{C}(\frac{R}{2},R))}\\
&\leq &C \| \vn\theta_R \|_{L^{r(\cdot)}(\mathcal{C}(\frac{R}{2},R))}\|P \|_{L^{\frac{p(\cdot)}{2}}(\Rt)}\|\vu \|_{L^{p(\cdot)}(\Rt)}.
\end{eqnarray*}
Since we have by hypothesis that $\vu\in L^{p(\cdot)}(\Rt)$ and $P\in L^{\frac{p(\cdot)}{2}}(\Rt)$, we only need to study the quantity $\|\vn\theta_R\|_{L^{r(\cdot)}(\mathcal{C}(\frac{R}{2},R))}$ where the exponent $r(\cdot)$ satisfies (\ref{Definition_rConjugueTheo1}). Following the same arguments as above we obtain $\underset{R\to +\infty}{\lim}\|\vn \theta_R \|_{L^{r(\cdot)}(\mathcal{C}(\frac{R}{2},R))}=0$, and thus we have the limit $\displaystyle \lim_{R\to +\infty}\beta_2(R)=0$.\\
\end{itemize}
With these two limits at hand for the quantities $\beta_1(R)$ and $\beta_2(R)$, we easily deduce from (\ref{aug_31T1}) that $\displaystyle \lim_{R\to +\infty}\beta(R)=0$. \\
\end{itemize}
We have thus proven that the terms $\alpha(R)$ and $\beta(R)$ given in (\ref{ineq_Base}) tend to $0$ as $R\to+\infty$: the proof of Theorem \ref{Theoreme_1} is now complete.\hfill $\blacksquare$\\
\subsection{Proof of the Theorem \ref{Theoreme_2}} 
Following the main ideas of the Subsection \ref{SubsecGeneral}, we only need to prove the limits (\ref{LimitsAlphaBeta}) where the quantities $\alpha(R)$ and $\beta(R)$ are defined in (\ref{ineq_Base}).
\begin{itemize}
\item[{\bf 1)}]{\bf Control for $\alpha(R)$}.  Proceeding as in (\ref{18aveq32_T1}) we have the estimate
\begin{equation}\label{EstimationAlphaTheo2}
\alpha(R)=\int_{B_R} \Delta \theta_R \frac{|\vu|^2}{2} dx \leq C\|\Delta \theta_R \|_{L^{\mathfrak{q}(\cdot)}(\Rt)}\|\vu\|_{L^{\mathfrak{p}(\cdot)}(\Rt)}^2,
\end{equation}
where we applied the H\"older inequalities with $\frac{2}{\mathfrak{p}(\cdot)}+\frac{1}{\mathfrak{q}(\cdot)}=1$ (\emph{i.e.} $\mathfrak{q}(\cdot)=\tfrac{\mathfrak{p}(\cdot)}{\mathfrak{p}(\cdot)-2}$). Since the variable exponent $\mathfrak{p}(\cdot)$ satisfies the conditions (\ref{Def_VariableExpTheo2}), we thus have the following conditions for $\mathfrak{q}(\cdot)$:
\begin{equation}\label{Definition_qConjugueTheo2}
\mathfrak{q}(x)=
\begin{cases}
\frac{9}{5}<  \mathfrak{q}_{(\Rt\setminus \mathscr{S})}^- \leq  \mathfrak{q}_{(\Rt\setminus \mathscr{S})}(x)\leq  \mathfrak{q}_{(\Rt\setminus \mathscr{S})}^+< 3,\\[5mm]
\frac{6\gamma+3}{2\gamma+3}< \mathfrak{q}_{\mathscr{S}}^- \leq  \mathfrak{q}_{ \mathscr{S}}(x)\leq  \mathfrak{q}_{ \mathscr{S}}^+<\frac{9}{5},
\end{cases}
\end{equation}
where $\mathscr{S}=\{(x_1, x_2, x_3)\in \Rt: x_2^2+x_3^2\leq x_1^\gamma, \; x_1>0\}$ and $0<\gamma<1$.\\

Since we have by hypothesis that $\vu\in L^{\mathfrak{p}(\cdot)}(\Rt)$, as explained previously, we can focus ourselves on the quantity $\|\Delta \theta_R \|_{L^{\mathfrak{q}(\cdot)}(\Rt)}$. Again, by the support properties of this localizing function $\theta_R$ we can write (see (\ref{ConditionSupportTheta}) above):
$$\|\Delta\theta_R\|_{L^{\mathfrak{q}(\cdot)}(\Rt)}=\|\Delta\theta_R\|_{L^{\mathfrak{q}(\cdot)}(\mathcal{C}(\frac{R}{2},R))}.$$
Now, we denote by $\mathscr{S}_1$ and $\mathscr{S}_2$ the sets defined by 
$$\mathscr{S}_1=\mathcal{C}\left(\tfrac{R}{2}, R\right)\cap \mathscr{S}\qquad \mbox{and}\qquad \mathscr{S}_2=\mathcal{C}\left(\tfrac{R}{2}, R\right) \setminus   \mathscr{S},$$
and we obtain
\begin{equation}\label{EstimationAlphaThetaTheo2}
\| \Delta \theta_R \|_{L^{\mathfrak{q}(\cdot)}(\mathcal{C}(\frac{R}{2},R))}\leq\|\Delta \theta_R \|_{L^{\mathfrak{q}_{\mathscr{S}_1}(\cdot)}(\mathscr{S}_1)}+\| \Delta \theta_R \|_{L^{\mathfrak{q}_{\mathscr{S}_2}(\cdot)}(\mathscr{S}_2)}.
\end{equation}
For the first term of the right-hand side above, following the same ideas given in 
(\ref{EstimationNormeLinfiniTheta_R1})-(\ref{EstimationC01}), we have
\begin{eqnarray*}
\|\Delta \theta_R \|_{L^{\mathfrak{q}_{\mathscr{S}_1}(\cdot)}(\mathscr{S}_1)}&\leq &C\|\Delta \theta_R  \|_{L^{\infty}  (\mathscr{S}_1)}\|1\|_{L^{\mathfrak{q}_{\mathscr{S}_1}(\cdot)}(\mathscr{S}_1)}\leq C R^{-2}\|1\|_{L^{\mathfrak{q}_{\mathscr{S}_1}(\cdot)}(\mathscr{S}_1)}\\
&\leq&C R^{-2}\max \{|\mathcal{C}(\tfrac{R}{2}, R)\cap \mathscr{S}|^{\frac 1 {\mathfrak{q}_{\mathscr{S}}^- }}, |\mathcal{C}(\tfrac{R}{2}, R)\cap \mathscr{S}|^{\frac 1 {\mathfrak{q}_{\mathscr{S}}^+}}\}.
\end{eqnarray*}
Since we have 
$$\mathcal{C}(\tfrac{R}{2}, R)\cap \mathscr{S}\subset \mathscr{B}=\{(x_1, x_2, x_3)\in \Rt: x_2^2+x_3^2\leq x_1^\gamma, \; 0<x_1<R\},$$
by computing the volume of the solid of revolution $\mathscr{B}$ defined above, we have
$$|\mathcal{C}(\tfrac{R}{2}, R)\cap \mathscr{S}|\leq |\mathscr{B}|=C\int_{0}^R x_1^{2\gamma} dx_1=CR^{2\gamma+1},$$ 
and we can write
\begin{equation}\label{EstimationS1}
\|\Delta \theta_R \|_{L^{\mathfrak{q}_{\mathscr{S}_1}(\cdot)}(\mathscr{S}_1)}\leq C\max \{R^{-2+\frac {2\gamma+1}{\mathfrak{q}_{\mathscr{S}}^-}}, R^{-2+\frac {2\gamma+1}{\mathfrak{q}_{\mathscr{S}}^+}}\}.
\end{equation}
But since by (\ref{Definition_qConjugueTheo2}), we have $\frac{5}{9}<\frac 1 {\mathfrak{q}_{\mathscr{S}}^+}\leq \frac 1 {\mathfrak{q}_{\mathscr{S}}^-}<\frac{2\gamma+3}{6\gamma+3}$, it comes (recall that $0<\gamma<1$)
$$-2+\frac{2\gamma+1}{\mathfrak{q}_{\mathscr{S}}^+}\leq -2+\frac{2\gamma+1}{\mathfrak{q}_{\mathscr{S}}^-}<-2+(2\gamma+1)\frac{2\gamma+3}{6\gamma+3}<0,$$ 
and thus all the powers of $R>1$ in the formula (\ref{EstimationS1}) are negative, from which we easily deduce that 
$$\|\Delta \theta_R \|_{L^{\mathfrak{q}_{\mathscr{S}_1}(\cdot)}(\mathscr{S}_1)}\underset{R\to +\infty}{\longrightarrow}0.$$
For the second term in the right-hand side of (\ref{EstimationAlphaThetaTheo2}), since we have $\mathscr{S}_2=\mathcal{C}(\frac{R}{2}, R)\setminus \mathscr{S}\subset \Rt\setminus \mathscr{S}$, 
by the same arguments given in (\ref{EstimationNormeLinfiniTheta_R2})-(\ref{EstimationC2}) we obtain
$$\|\Delta \theta_R \|_{L^{\mathfrak{q}_{\mathscr{S}_2}(\cdot)}(\mathscr{S}_2)}\underset{R\to +\infty}{\longrightarrow}0.$$
With all these estimate at hand and coming back to (\ref{EstimationAlphaTheo2}) we finally obtain that 
$$\underset{R\to+\infty}{\lim}\alpha(R)\to 0.$$

\item[{\bf 2)}]{\bf Control for $\beta(R)$}. Following the ideas in (\ref{aug_31T1}) we have
\begin{equation}\label{EstimationGeneralBeta12Theo2}
\left|\beta(R)\right| \leq \frac{1}{2} \underbrace{\int_{\mathcal{C}\left(\frac{R}{2}, R\right)}|\vn \theta_R||\vu|^3dx}_{\beta_1(R)}+\underbrace{\int_{\mathcal{C}\left(\frac{R}{2}, R\right)}|\vn \theta_R||P \| \vu| dx}_{\beta_2(R)}.
\end{equation}
For the first term in the right-hand side above, following (\ref{2mayeq3T1}) we write:
\begin{equation}\label{EstimationBeta1Theo2}
\beta_1(R)\leq C\|\vn\theta_R \|_{L^{\mathfrak{r}(\cdot)}(\mathcal{C}\left(\frac{R}{2}, R\right))}\|\vu \|_{L^{\mathfrak{p}(\cdot)}(\Rt)}^3,
\end{equation}
where we used the H\"older inequality with $\frac{3}{\mathfrak{p}(\cdot)}+\frac{1}{\mathfrak{r}(\cdot)}=1$ (\emph{i.e.} we have $\mathfrak{r}(\cdot)=\frac{\mathfrak{p}(\cdot)}{\mathfrak{p}(\cdot)-3}$). Since the exponent $\mathfrak{p}(\cdot)$ is now driven by the conditions stated in (\ref{Def_VariableExpTheo2}), the variable exponent $\mathfrak{r}(\cdot)$ satisfies: 
\begin{equation}\label{Definition_rConjugueTheo2}
\mathfrak{r}(x)=
\begin{cases}
3<  \mathfrak{r}_{(\Rt\setminus \mathscr{S})}^- \leq  \mathfrak{r}_{(\Rt\setminus \mathscr{S})}(x)\leq  \mathfrak{r}_{(\Rt\setminus \mathscr{S})}^+< +\infty,\\[5mm]
2\gamma+1< \mathfrak{r}_{\mathscr{S}}^- \leq  \mathfrak{r}_{ \mathscr{S}}(x)\leq  \mathfrak{r}_{ \mathscr{S}}^+<3.
\end{cases}
\end{equation}
Again, since by hypothesis we have $\|\vu \|_{L^{\mathfrak{p}(\cdot)}(\Rt)}<+\infty$, to understand the behavior of $\beta_1(R)$ we only need to study the quantity $\|\vn\theta_R \|_{L^{\mathfrak{r}(\cdot)}(\mathcal{C}\left(\frac{R}{2}, R\right))}$. Following the main ideas used to obtain the estimate (\ref{EstimationBeta1C}) and by the same arguments displayed above (recall that $|\mathcal{C}(\tfrac{R}{2}, R)\cap \mathscr{S}|\leq CR^{2\gamma+1}$), we have
$$\| \vn \theta_R \|_{L^{\mathfrak{r}(\cdot)}(\mathcal{C}(\frac{R}{2},R))}\leq C\max \{R^{-1+\frac{2\gamma+1}{\mathfrak{r}_{\mathscr{S}}^-}}, R^{-1+\frac{2\gamma+1}{\mathfrak{r}_{\mathscr{S}}^+}}\}+C\max\{R^{-1+\frac3{\mathfrak{r}_{(\Rt\setminus\mathscr{S})}^-}},R^{-1+\frac 3{\mathfrak{r}_{(\Rt\setminus\mathscr{S})}^+}}\}.$$
Now, by the restrictions (\ref{Definition_rConjugueTheo2}) we easily deduce that 
$$-1+\frac{2\gamma+1}{\mathfrak{r}_{\mathscr{S}}^+}\leq -1+\frac{2\gamma+1}{\mathfrak{r}_{\mathscr{S}}^-}<0\qquad \mbox{and}\qquad -1+\frac 3{\mathfrak{r}_{(\Rt\setminus\mathscr{S})}^+}\leq -1+\frac 3{\mathfrak{r}_{(\Rt\setminus\mathscr{S})}^-}<0,$$
and since all the powers of the parameter $R>1$ are negative we have $\|\vn \theta_R \|_{L^{\mathfrak{r}(\cdot)}(\mathcal{C}(\frac{R}{2},R))}\underset{R\to+\infty}{\longrightarrow}0$. Thus, getting back to (\ref{EstimationBeta1Theo2}) we have $\displaystyle \lim_{R\to +\infty}\beta_1(R)=0$.\\

For the term $\beta_2(R)$ given in (\ref{EstimationGeneralBeta12Theo2}) we proceed as follows: by the H\"older inequalities with $\frac{1}{\mathfrak{p}(\cdot)}+\frac{2}{\mathfrak{p}(\cdot)}+\frac{1}{\mathfrak{r}(\cdot)}=1$, we obtain
\begin{eqnarray*}
\beta_2(R)&=&\int_{\mathcal{C}\left(\frac{R}{2}, R\right)}|\vn \theta_R||P \| \vu| dx\leq C \| \vn\theta_R \|_{L^{\mathfrak{r}(\cdot)}(\mathcal{C}(\frac{R}{2},R))} \|P \|_{L^{\frac{\mathfrak{p}(\cdot)}{2}}(\mathcal{C}(\frac{R}{2},R))}\|\vu \|_{L^{\mathfrak{p}(\cdot)}(\mathcal{C}(\frac{R}{2},R))}\\
&\leq &C \| \vn\theta_R \|_{L^{\mathfrak{r}(\cdot)}(\mathcal{C}(\frac{R}{2},R))}\|P \|_{L^{\frac{\mathfrak{p}(\cdot)}{2}}(\Rt)}\|\vu \|_{L^{\mathfrak{p}(\cdot)}(\Rt)}.
\end{eqnarray*}
Now, since we have by hypothesis that $\vu\in L^{\mathfrak{p}(\cdot)}(\Rt)$ and $P\in L^{\frac{\mathfrak{p}(\cdot)}{2}}(\Rt)$, as before, we only need to study the quantity $\|\vn\theta_R\|_{L^{\mathfrak{r}(\cdot)}(\mathcal{C}(\frac{R}{2},R))}$ where the exponent $\mathfrak{r}(\cdot)$ satisfies this time the relationships $\frac{1}{\mathfrak{p}(\cdot)}+\frac{2}{\mathfrak{p}(\cdot)}+\frac{1}{\mathfrak{r}(\cdot)}=1$, \emph{i.e.} $\mathfrak{r}(\cdot)=\frac{\mathfrak{p}(\cdot)}{\mathfrak{p}(\cdot)-3}$, where $\mathfrak{r}(\cdot)$ verifies (\ref{Definition_rConjugueTheo2}). Thus, following exaclty the same arguments as above we obtain $\underset{R\to +\infty}{\lim}\|\vn \theta_R \|_{L^{\mathfrak{r}(\cdot)}(\mathcal{C}(\frac{R}{2},R))}=0$, and we finally obtain the limit $\displaystyle\lim_{R\to +\infty}\beta_2(R)=0$.\\

We have obtained that $\displaystyle\lim_{R\to +\infty}\beta_1(R)=\displaystyle\lim_{R\to +\infty}\beta_2(R)=0$, from which, getting back to the expression (\ref{EstimationGeneralBeta12Theo2}), we easily deduce that $\displaystyle\lim_{R\to +\infty}\beta(R)=0$.\\
\end{itemize}

After the previous estimates for the terms $\alpha(R)$ and $\beta(R)$, we can consider the control (\ref{ineq_Base}), and we deduce, by (\ref{LimiteSobolev}), that $\|\vu\|_{\dot{H}^1}=0$ and that $\vu\equiv 0$. The proof of the Theorem \ref{Theoreme_2} is now finished. \hfill $\blacksquare$

\subsection{Proof of the Theorem \ref{Theoreme_3}} 
As in the proof of the two previous results, our starting point is the inequality (\ref{ineq_Base}) given by
$$\int_{B_{\frac{R}{2}}}|\vn \otimes \vu|^2 dx \leq \underbrace{\int_{B_R} \Delta \theta_R \frac{|\vu|^2}{2} dx}_{\alpha(R)}+\underbrace{\int_{B_R} \vn \theta_R \cdot\left(\left(\frac{|\vu|^2}{2}+P\right) \vu\right)dx}_{\beta(R)},$$
and we aim to prove that $\displaystyle\lim_{R\to +\infty}\alpha(R)=\displaystyle\lim_{R\to +\infty}\beta(R)=0$.
\begin{itemize}
\item[{\bf 1)}]{\bf Control for $\alpha(R)$}.  By the H\"older inequalities with\footnote{Recall that the H\"older inequalities are still fully available in this context with the values $1\leq {\bf p}(\cdot)\leq +\infty$, see the book \cite[Theorem 2.26]{Cruz_Libro} for more details.}  $\frac{2}{{\bf p}(\cdot)}+\frac{1}{{\bf q}(\cdot)}=1$, we obtain 
$$\alpha(R)=\int_{B_R} \Delta \theta_R \frac{|\vu|^2}{2} dx \leq C\|\Delta \theta_R \|_{L^{{\bf q}(\cdot)}(\Rt)}\|\vu\|_{L^{{\bf p}(\cdot)}(\Rt)}^2.$$
As before, since we have by hypothesis that $\|\vu\|_{L^{{\bf p}(\cdot)}(\Rt)}<+\infty$, in order to study the behavior as $R\to+\infty$ of the quantity $\alpha(R)$, we only need to study the term $\|\Delta \theta_R \|_{L^{{\bf q}(\cdot)}(\Rt)}$.\\

Now, as the variable exponent ${\bf p}(\cdot)$ satisfies the conditions (\ref{Def_VariableExpTheo3}), we deduce that the variable exponent ${\bf q}(\cdot)=\tfrac{{\bf p}(\cdot)}{{\bf p}(\cdot)-2}$ is given by 
\begin{equation}\label{Definition_qConjugueTheo3}
{\bf q}(x)=
\begin{cases}
\frac{9}{5}<  {\bf q}_{(\Rt\setminus \mathscr{N})}^- \leq  {\bf q}_{(\Rt\setminus \mathscr{N})}(x)\leq  {\bf q}_{(\Rt\setminus \mathscr{N})}^+< 3,\\[5mm]
{\bf q}_{ \mathscr{N}}(x)={\bf q}_{ \mathscr{N}}^-={\bf q}_{ \mathscr{N}}^+=1,
\end{cases}
\end{equation}
where $\mathscr{N}=\{(x_1, x_2, x_3)\in \Rt: x_2^2+x_3^2\leq x_1^{-\frac{\sigma}{2}},\;\; x_1>0\}$ and $0<\sigma<1$. We thus write, 
\begin{equation}\label{EtudeLaplacienTheta3}
\|\Delta \theta_R \|_{L^{{\bf q}(\cdot)}(\Rt)}=\|\Delta\theta_R\|_{L^{{\bf q}(\cdot)}(\mathcal{C}(\frac{R}{2},R))}\leq\|\Delta \theta_R \|_{L^{{\bf q}_{\mathscr{N}_1}(\cdot)}(\mathscr{N}_1)}+\| \Delta \theta_R \|_{L^{{\bf q}_{\mathscr{N}_2}(\cdot)}(\mathscr{N}_2)}.
\end{equation}
with $\mathscr{N}_1=\mathcal{C}\left(\tfrac{R}{2}, R\right)\cap \mathscr{N}$ and $\mathscr{N}_2=\mathcal{C}\left(\tfrac{R}{2}, R\right) \setminus \mathscr{N}$.\\

For the term $\|\Delta \theta_R \|_{L^{{\bf q}_{\mathscr{N}_1}(\cdot)}(\mathscr{N}_1)}$ given in (\ref{EtudeLaplacienTheta3}), we write
$$\|\Delta \theta_R \|_{L^{{\bf q}_{\mathscr{N}_1}(\cdot)}(\mathscr{N}_1)} \leq C R^{-2}\|1\|_{L^{{\bf q}_{\mathscr{N}_1}(\cdot)}(\mathscr{N}_1)}\leq C R^{-2}\max \{|\mathcal{C}(\tfrac{R}{2}, R)\cap \mathscr{N}|^{\frac 1 {{\bf q}_{\mathscr{N}}^- }}, |\mathcal{C}(\tfrac{R}{2}, R)\cap \mathscr{N}|^{\frac 1 {{\bf q}_{\mathscr{N}}^+}}\}.$$
Now, since $\mathcal{C}(\tfrac{R}{2}, R)\cap \mathscr{N}\subset \mathscr{N}$  we thus have by (\ref{Definition_qConjugueTheo3}) that ${\bf q}_{\mathscr{N}}^- ={\bf q}_{\mathscr{N}}^- =1$, moreover we have the set inclusion
$$\mathcal{C}(\tfrac{R}{2}, R)\cap \mathscr{N}\subset \mathscr{D}=\{(x_1, x_2, x_3)\in \Rt: x_2^2+x_3^2\leq x_1^{-\frac{\sigma}{2}},\;\;0<x_1<R\},$$
and computing the volume of the solid of revolution of the set $\mathscr{D}$ we obtain the following estimate (recall that $0<\sigma<1$):
$$|\mathcal{C}(\tfrac{R}{2}, R)\cap \mathscr{N}|\leq |\mathscr{D}|=C\int_{0}^{R}x_1^{-\sigma}dx_1= CR^{1-\sigma},$$ 
and we can write
$$\|\Delta \theta_R \|_{L^{{\bf q}_{\mathscr{N}_1}(\cdot)}(\mathscr{N}_1)} \leq CR^{-1-\sigma}\underset{R\to+\infty}{\longrightarrow}0.$$
We study now the term $\| \Delta \theta_R \|_{L^{{\bf q}_{\mathscr{N}_2}(\cdot)}(\mathscr{N}_2)}$  given in (\ref{EtudeLaplacienTheta3}). Since we have the set inclusions $\mathscr{N}_2=\mathcal{C}(\frac{R}{2}, R)\setminus \mathscr{N}\subset \Rt\setminus \mathscr{N}$ and since over this set we have that the variable exponent ${\bf q}(\cdot)$ satisfies the condition (\ref{Definition_qConjugueTheo3}), by the same arguments given in (\ref{EstimationNormeLinfiniTheta_R2})-(\ref{EstimationC2}) we obtain
$$\|\Delta \theta_R \|_{L^{{\bf q}_{\mathscr{N}_2}(\cdot)}(\mathscr{N}_2)}\underset{R\to +\infty}{\longrightarrow}0.$$
With these two estimates at hand we deduce, coming back to (\ref{EtudeLaplacienTheta3}), that $\|\Delta \theta_R \|_{L^{{\bf q}(\cdot)}(\Rt)}\underset{R\to +\infty}{\longrightarrow}0$ and this fact implies that we have
$$\alpha(R)\underset{R\to +\infty}{\longrightarrow}0.$$

\item[{\bf 2)}]{\bf Control for $\beta(R)$}. The control of this terms follows the same general ideas used in (\ref{aug_31T1}), we can thus write
\begin{equation}\label{EstimationGeneralBeta12Theo3}
\left|\beta(R)\right| \leq \frac{1}{2} \underbrace{\int_{\mathcal{C}\left(\frac{R}{2}, R\right)}|\vn \theta_R||\vu|^3dx}_{\beta_1(R)}+\underbrace{\int_{\mathcal{C}\left(\frac{R}{2}, R\right)}|\vn \theta_R||P \| \vu| dx}_{\beta_2(R)}.
\end{equation}
For the first term in the right-hand side above, following (\ref{2mayeq3T1}) we write:
\begin{equation}\label{EstimationBeta1Theo3}
\beta_1(R)\leq C\|\vn\theta_R \|_{L^{{\bf r}(\cdot)}(\mathcal{C}\left(\frac{R}{2}, R\right))}\|\vu \|_{L^{{\bf p}(\cdot)}(\Rt)}^3,
\end{equation}
where we used the H\"older inequality with $\frac{3}{{\bf p}(\cdot)}+\frac{1}{{\bf r}(\cdot)}=1$. Since the exponent ${\bf p}(\cdot)$ is given by the conditions stated in (\ref{Def_VariableExpTheo3}), the variable exponent ${\bf r}(\cdot)=\frac{{\bf p}(\cdot)}{{\bf p}(\cdot)-3}$ satisfies: 
\begin{equation}\label{Definition_rConjugueTheo3}
{\bf r}(x)=
\begin{cases}
3<  {\bf r}_{(\Rt\setminus \mathscr{N})}^- \leq  {\bf r}_{(\Rt\setminus \mathscr{N})}(x)\leq  {\bf r}_{(\Rt\setminus \mathscr{N})}^+< +\infty,\\[5mm]
{\bf r}_{\mathscr{N}}^- =  {\bf r}_{ \mathscr{N}}(x)=  {\bf r}_{ \mathscr{N}}^+=1.
\end{cases}
\end{equation}
 To understand the behavior of $\beta_1(R)$, since by hypothesis we have $\|\vu \|_{L^{{\bf p}(\cdot)}(\Rt)}<+\infty$, we only need to study the quantity $\|\vn\theta_R \|_{L^{{\bf r}(\cdot)}(\mathcal{C}\left(\frac{R}{2}, R\right))}$. Following the main ideas used to obtain the estimate (\ref{EstimationBeta1C}) we have
\begin{eqnarray*}
\| \vn \theta_R \|_{L^{{\bf r}(\cdot)}(\mathcal{C}(\frac{R}{2},R))}&\leq &CR^{-1}\max \{|\mathcal{C}(\tfrac{R}{2}, R)\cap \mathscr{N}|^{\frac{1}{{\bf r}_{\mathscr{N}}^-}},|\mathcal{C}(\tfrac{R}{2}, R)\cap \mathscr{N}|^{\frac{1}{{\bf r}_{\mathscr{N}}^+}}\}\\
&&+CR^{-1}\max\{|\mathcal{C}(\tfrac{R}{2}, R)\setminus \mathscr{N}|^{\frac1{{\bf r}_{(\Rt\setminus\mathscr{N})}^-}},|\mathcal{C}(\tfrac{R}{2}, R)\setminus \mathscr{N}|^{\frac 1{{\bf r}_{(\Rt\setminus\mathscr{N})}^+}}\}.
\end{eqnarray*}
Recalling that we have $|\mathcal{C}(\tfrac{R}{2}, R)\cap \mathscr{N}|\leq CR^{1-\sigma}$ and $|\mathcal{C}(\tfrac{R}{2}, R)\setminus \mathscr{N}|\leq CR^3$, we obtain, since ${\bf r}_{\mathscr{N}}^- =  {\bf r}_{ \mathscr{N}}(x)=  {\bf r}_{ \mathscr{N}}^+=1$: 
\begin{eqnarray*}
\| \vn \theta_R \|_{L^{{\bf r}(\cdot)}(\mathcal{C}(\frac{R}{2},R))}&\leq &CR^{-\sigma}+C\max\{R^{-1+\frac3{{\bf r}_{(\Rt\setminus\mathscr{N})}^-}},R^{-1+\frac 3{{\bf r}_{(\Rt\setminus\mathscr{N})}^+}}\}
\end{eqnarray*}
Now, by the restrictions (\ref{Definition_rConjugueTheo3}) we easily deduce that 
$ -1+\frac 3{{\bf r}_{(\Rt\setminus\mathscr{N})}^+}\leq -1+\frac 3{{\bf r}_{(\Rt\setminus\mathscr{N})}^-}<0$ and since all the powers of the parameter $R>1$ are negative we have $\|\vn \theta_R \|_{L^{{\bf r}(\cdot)}(\mathcal{C}(\frac{R}{2},R))}\underset{R\to+\infty}{\longrightarrow}0$. Thus, getting back to (\ref{EstimationBeta1Theo3}) we have $\displaystyle \lim_{R\to +\infty}\beta_1(R)=0$.\\

For the term $\beta_2(R)$ given in (\ref{EstimationGeneralBeta12Theo3}) we proceed following the same ideas used previously: by the H\"older inequalities with $\frac{1}{{\bf p}(\cdot)}+\frac{2}{{\bf p}(\cdot)}+\frac{1}{{\bf r}(\cdot)}=1$, we obtain
\begin{eqnarray*}
\beta_2(R)&=&\int_{\mathcal{C}\left(\frac{R}{2}, R\right)}|\vn \theta_R||P \| \vu| dx\leq C \| \vn\theta_R \|_{L^{{\bf r}(\cdot)}(\mathcal{C}(\frac{R}{2},R))} \|P \|_{L^{\frac{{\bf p}(\cdot)}{2}}(\mathcal{C}(\frac{R}{2},R))}\|\vu \|_{L^{{\bf p}(\cdot)}(\mathcal{C}(\frac{R}{2},R))}\\
&\leq &C \| \vn\theta_R \|_{L^{{\bf r}(\cdot)}(\mathcal{C}(\frac{R}{2},R))}\|P \|_{L^{\frac{{\bf p}(\cdot)}{2}}(\Rt)}\|\vu \|_{L^{{\bf p}(\cdot)}(\Rt)}.
\end{eqnarray*}
Since by hypothesis $\vu\in L^{{\bf p}(\cdot)}(\Rt)$ and $P\in L^{\frac{{\bf p}(\cdot)}{2}}(\Rt)$, as before, we only need to study the quantity $\|\vn\theta_R\|_{L^{{\bf r}(\cdot)}(\mathcal{C}(\frac{R}{2},R))}$ where the variable exponent ${\bf r}(\cdot)=\frac{{\bf p}(\cdot)}{{\bf p}(\cdot)-3}$ is given by (\ref{Definition_rConjugueTheo3}). Thus, following exaclty the same arguments as above we obtain $\underset{R\to +\infty}{\lim}\|\vn \theta_R \|_{L^{{\bf r}(\cdot)}(\mathcal{C}(\frac{R}{2},R))}=0$, and we finally obtain the limit $\displaystyle\lim_{R\to +\infty}\beta_2(R)=0$.\\

We have obtained that $\displaystyle\lim_{R\to +\infty}\beta_1(R)=\displaystyle\lim_{R\to +\infty}\beta_2(R)=0$, from which, getting back to the expression (\ref{EstimationGeneralBeta12Theo3}), we easily deduce that $\displaystyle\lim_{R\to +\infty}\beta(R)=0$.\\
\end{itemize}
We have thus proven that the terms $\alpha(R)$ and $\beta(R)$ given in (\ref{ineq_Base}) tend to $0$ as $R\to+\infty$: the proof of Theorem \ref{Theoreme_3} is now complete.\hfill $\blacksquare$\\

\paragraph{\bf Conflict of interest.} On behalf of all authors, the corresponding author states that there is no conflict of interest.
\paragraph{\bf Acknowledgements.} 
The second author is supported by the ANID postdoctoral program BCH 2022 grant No. 74220003.
\paragraph{\bf Contribution.} 
All the authors have equally contributed to this research.  


\end{document}